\documentclass[11pt,a4paper,leqno,twoside]{amsart}
\usepackage{latexsym,amssymb,amsmath}
\input xy
\xyoption{all}

\def\CC{\mathbb C}

\def\RR{\mathbb R}
\def\HH{\mathbb H}
\def\AA{{\mathbb A}}

\def\OO{\mathbb O}

\def\11{\mathbf 1}
\def\PP{\mathbb P}

\def\e1{\varepsilon_1}
\def\e2{\varepsilon_2}
\def\e3{\varepsilon_3}

\def\P2{{\PP}^2}

\def\00{\underline{0}}
\def\J0{{\cal J}_3(\underline{0})}

\def\PJ0{\PP({\cal J}_3(\underline{0}))}

\def\e{\varepsilon}

\def\AP2{{\AA\PP}^2}
\def\RP2{{\RR\PP}^2}
\def\CP2{{\CC\PP}^2}
\def\HP2{{\HH\PP}^2}
\def\OP2{{\OO\PP}^2}

\newtheorem{theo}{Theorem}[section]
\newtheorem{coro}[theo]{Corollary}
\newtheorem{lemm}[theo]{Lemma}
\newtheorem{prop}[theo]{Proposition}

\theoremstyle{definition}

\theoremstyle{remark}
\newtheorem{rema}[theo]{Remark}

\begin{document}
\title[Henselian discrete fields]{Henselian discrete valued stable fields}

\author{I.D. Chipchakov}
\address{Institute of Mathematics and Informatics\\Bulgarian Academy
of Sciences\\1113 Sofia, Bulgaria: E-mail address:
chipchak@math.bas.bg}

\dedicatory{Dedicated to Vesselin Drensky on the occasion of his 70th
anniversary}

\begin{abstract}
Let $(K, v)$ be a Henselian discrete valued field with residue field
$\widehat K$ of characteristic $q \ge 0$, and Brd$_{p}(K)$ be the Brauer
$p$-dimension of $K$, for each prime $p$. The present paper
shows that if $p = q$, then Brd$_{p}(K) \le 1$ if and only if
$\widehat K$ is a $p$-quasilocal field and the degree $[\widehat
K\colon \widehat K ^{p}]$ is $\le p$. This complements our earlier
result that, in case $p \neq q$, we have Brd$_{p}(K) \le 1$ if and
only if $\widehat K$ is $p$-quasilocal and Brd$_{p}(\widehat K) \le
1$.
\end{abstract}

\maketitle

\section{Introduction}
\par
\medskip
Let $E$ be a field, $E _{\rm sep}$ its separable closure, Br$(E)$
the Brauer group of $E$, $s(E)$ the class of associative
finite-dimensional central simple algebras over $E$, and $d(E)$
the subclass of division algebras $D \in s(E)$. For each $A \in
s(E)$, let deg$(A)$, ind$(A)$ and exp$(A)$ be the degree, the
Schur index and the exponent of $A$, respectively. It is
well-known (cf. \cite{P}, Sect. 14.4) that exp$(A)$ divides
ind$(A)$ and shares with it the same set of prime divisors; also,
ind$(A) \mid {\rm deg}(A)$, and deg$(A) = {\rm ind}(A)$ if and
only if $A \in d(E)$. Note that ind$(B _{1} \otimes _{E} B _{2}) =
{\rm ind}(B _{1}){\rm ind}(B _{2})$ if $B _{1}, B _{2} \in s(E)$
and g.c.d.$\{{\rm ind}(B _{1}), {\rm ind}(B _{2})\} = 1$;
equivalently, $B _{1} ^{\prime } \otimes _{E} B _{2} ^{\prime }
\in d(E)$ in case $B _{j} ^{\prime } \in d(E)$, $j = 1, 2$, and
g.c.d.$\{{\rm deg}(B _{1} ^{\prime }), {\rm deg}(B _{2} ^{\prime
})\}$ $= 1$ (see \cite{P}, Sect. 13.4). Since Br$(E)$ is an
abelian torsion group and ind$(A)$, exp$(A)$ are invariants both
of $A$ and its equivalence class $[A] \in {\rm Br}(E)$, these
results reduce the study of the restrictions on the pairs
ind$(A)$, exp$(A)$, $A \in s(E)$, to the special case of
$p$-primary pairs, for an arbitrary fixed prime $p$. The Brauer
$p$-dimensions Brd$_{p}(E)$, $p \in \mathbb P$, where $\mathbb P$
is the set of prime numbers, contain essential information on
these restrictions. We say that Brd$_{p}(E) = n < \infty $, for a
given $p \in \mathbb P$, if $n$ is the least integer $\ge 0$, for
which ind$(A _{p}) \mid {\rm exp}(A _{p}) ^{n}$ whenever $A _{p}
\in s(E)$ and $[A _{p}]$ lies in the $p$-component Br$(E) _{p}$ of
Br$(E)$; if no such $n$ exists, we put Brd$_{p}(E) = \infty $. For
instance, Brd$_{p}(E) \le 1$, for all $p \in \mathbb P$, if and
only if $E$ is a stable field, i.e. deg$(D) = {\rm exp}(D)$, for
each $D \in d(E)$; Brd$_{p'}(E) = 0$, for some $p ^{\prime } \in
\mathbb P$, if and only if Br$(E) _{p'} = \{0\}$. The absolute
Brauer $p$-dimension abrd$_{p}(E)$ of $E$ is defined to be the
supremum of Brd$_{p}(R)\colon R \in {\rm Fe}(E)$, where Fe$(E)$ is
the set of finite extensions of $E$ in $E _{\rm sep}$. We have
abrd$_{p}(E) \le 1$, $p \in \mathbb P$, if $E$ is an absolutely
stable field, i.e. its finite extensions are stable fields.
Important fields of this kind have been exhibited by class field
theory and the theory of algebraic surfaces, which show that
Brd$_{p}(\Phi ) = {\rm abrd}_{p}(\Phi ) = 1$, $p \in \mathbb P$,
if $\Phi $ is a global or local field (see, e.g., \cite{Re},
(31.4) and (32.19)), or a finitely-generated extension of
transcendence degree $2$ over an algebraically closed field $\Phi
_{0}$ \cite{Jon, Li}.
\par
Similarly to other topics in the theory of central simple algebras
and Brauer groups of fields (see, e.g., \cite{TW}, Chaps. 9-12),
the study of the sequence Brd$_{p}(E), {\rm abrd}_{p}(E)$, $p \in
\mathbb P$, brings useful general results if it restricts at one
point or another to certain classes of Henselian (valued) fields
or other suitably chosen special fields. The restriction on $E$
allows to find formulae for Brd$_{p}(E)$ and abrd$_{p}(E)$, and to
use them for constructing fields $E ^{\prime }$ with prescribed
sequences Brd$_{p}(E ^{\prime }), {\rm abrd}_{p}(E ^{\prime })$,
$p \in \mathbb P$ \cite{Ch7}. This in turn provides new
information on the behaviour of index-exponent relations under
finitely-generated field extensions (see, e.g., the answer to
\cite{ABGV}, Problem~4.4, given in \cite{Ch4, Ch5}, or else
\cite{Ch7}, Corollary~5.6 and \cite{Ch5}, Remark~5.5, for the
ground fields in a sequence of examples disproving Conjecture~2 of
\cite{PC}). The chosen approach also contributes to better
knowledge of Brauer groups of absolutely stable fields, viewed as
abstract abelian torsion groups \cite{Ch1}, Corollary~4.7. More
recently, it has been shown in \cite{Br} that some absolutely
stable fields (with absolute stability proved in \cite{Ch1}) admit
noncyclic division algebras of degree $2 ^{\nu }$, for every
integer $\nu \ge 2$.
\par
A nontrivial Krull valuation $v$ of a field $K$ is called
Henselian, if it extends uniquely, up-to equivalence, to a
valuation $v _{L}$ on each algebraic extension $L$ of $K$. The
stability condition on a Henselian (valued) field $(K _{0}, v
_{0})$ with a residue field $\widehat K _{0}$ of zero
characteristic has been fully characterized in \cite{Ch1} by
conditions on $\widehat K _{0}$ and the value group $v
_{0}(K_{0})$. Also, \cite{Ch6}, Proposition~3.5, and results of
\cite{Ch1} characterize maximally complete stable fields $(K _{q},
v _{q})$ with $\widehat K _{q}$ perfect and char$(K _{q}) = q >
0$. For example, by \cite{Ch1}, Corollary~4.5 (ii), the iterated
formal (Laurent) power series field $J((X))((Y))$ in $2$ variables
over a field $J$ is absolutely stable if and only if $J$ is
perfect and the absolute Galois group $\mathcal{G}_{J} :=
\mathcal{G}(J _{\rm sep}/J)$ is metabelian of cohomological
dimension cd$(\mathcal{G}_{J}) \le 1$, in the sense of \cite{S1}.
By Lemma~1.2 of \cite{Ch2}, $\mathcal{G}_{J}$ possesses the noted
properties if and only if its Sylow pro-$p$-groups are
topologically isomorphic to the additive group $\mathbb Z _{p}$ of
$p$-adic integers whenever $p \in \mathbb P$ and the cohomological
$p$-dimension cd$_{p}(\mathcal{G}_{J})$ (in the sense of
\cite{S1}) is nonzero. Therefore, $J((X))((Y))$ is absolutely
stable if the field $J$ is quasifinite, i.e. perfect with
$\mathcal{G}_{J}$ isomorphic to the topological group product
$\prod _{p \in \mathbb P} \mathbb Z _{p}$ (see also Remark
\ref{rema5.3} (ii)).
\par
\medskip The present paper can be viewed as a continuation of
\cite{Ch1}. It completes the characterization of Henselian
discrete valued (abbr, HDV) stable fields by properties of their
residue fields. Combined with \cite{PS}, Theorem~2, it determines
Brd$_{p}(K)$ in case $(K, v)$ is an HDV-field, char$(\widehat K) =
p > 0$, and the degree $[\widehat K\colon \widehat K ^{p}]$ is at
most equal to $p$, where $\widehat K ^{p} = \{\hat \alpha
^{p}\colon \hat \alpha \in \widehat K\}$.

\section{Statements of the main results}
\label{Sec:2} It is known that Brd$_{p}(\widehat K) \le {\rm
Brd}_{p}(K)$, $p \in \mathbb P$, for any Henselian field $(K, v)$
(see Theorem~2.8 of \cite{JW}, or Lemma \ref{lemm3.3} below).
Therefore, $\widehat K$ is a stable field, provided that so is
$K$. The problem of characterizing Henselian stable fields with
$p$-indivisible value groups is related to the study of
$p$-quasilocal fields. By a $p$-quasilocal field, for some $p \in
\mathbb P$, we mean a field $E$ satisfying one of the following
two conditions: Brd$_{p}(E) = 0$ or $E(p) = E$, where $E(p)$ is
the maximal $p$-extension of $E$ (in $E _{\rm sep}$); Brd$_{p}(E)
\neq 0$, $E(p) \neq E$, and every extension of $E$ in $E(p)$ of
degree $p$ is embeddable as an $E$-subalgebra in each $D _{p} \in
d(E)$ of degree $p$. We say that the field $E$ is quasilocal, if
its finite extensions are $p$-quasilocal fields, for every $p \in
\mathbb P$. The class of quasilocal fields contains local fields;
in addition, it is essentially larger than the class of HDV-fields
with quasifinite residue fields (cf. \cite{Se}, Ch. XIII, Sect. 3,
and \cite{Ch3}, Remark~3.7). As to global fields, they are not
$p$-quasilocal, for any $p \in \mathbb{P}$; if $F$ is a global
field, then one obtains from the Grunwald-Wang theorem and the
description of Br$(F)$ by class field theory (see \cite{AT}, Ch.
X, and \cite{We}, Ch. XIII, Sects. 3 and 6, respectively) that for
any $p \in \mathbb{P}$ and each $\Delta _{p} \in d(F)$ with
deg$(\Delta _{p}) = p$, there exist infinitely many extensions
$\Phi _{p}$ of $F$ in $F(p)$, such that $[\Phi _{p}\colon F] = p$
and $\Delta _{p} \otimes _{F} \Phi _{p} \in d(\Phi _{p})$. Note
also that quasilocal fields are absolutely stable \cite{Ch3},
Proposition~2.3. This is implied by the fact that a $p$-quasilocal
field $E$ satisfies Brd$_{p}(E) \le 1$ in the following cases: (i)
$E(p) \neq E$ \cite{Ch3}, Theorem~3.1 (ii); (ii) $E$ contains a
primitive $p$-th root of unity; (iii) char$(E) = p$. Moreover, in
cases (ii) and (iii), the assumption that Brd$_{p}(E) \neq 0$
ensures that $E(p) \neq E$ (see \cite{MS}, (16.2), and \cite{A1},
Ch. VII, Theorem~28). Henselian stable fields and $p$-quasilocal
fields are related by the following results:
\par
\medskip
\begin{prop}
\label{prop2.1} Let $(K, v)$ be a Henselian field and $p$ be a
prime. Then:
\par
{\rm (a)} $\widehat K$ is a $p$-quasilocal field if $v(K) \neq
pv(K)$ and {\rm Brd}$_{p}(K) \le 1$.
\par
{\rm (b)} {\rm Brd}$_{p}(K) \le 1$, provided that $\widehat K$ is
$p$-quasilocal, $p \neq {\rm char}(\widehat K)$, {\rm
Brd}$_{p}(\widehat K) \le 1$, and the quotient group $v(K)/pv(K)$
has order $p$.
\end{prop}
\par
\medskip\noindent
Proposition \ref{prop2.1} (a) follows from \cite{Ch3},
Proposition~2.1 (with its proof). Proposition \ref{prop2.1} (b) is
a special case of \cite{Ch6}, Theorem~4.1; it can also be deduced
from \cite{Ch1}, Theorem~3.1~(a). Using Proposition \ref{prop2.1}
(a), (b) and well-known results about the reduction of Schur
indices and exponents under a scalar extension of finite degree
over the centre (cf. \cite{P}, Sects. 13.4 and 14.4), one sees
that the problem of characterizing stable HDV-fields reduces to
the one of finding a necessary and sufficient condition that
Brd$_{p}(K) \le 1$, where $(K, v)$ is an HDV-field with
char$(\widehat K) = p$. The main result of \cite{Ch8}, stated
below, takes a step towards achieving this goal. It shows that if
Brd$_{p}(K) \le 1$, then $[\widehat K\colon \widehat K ^{p}] \le
p$ (in case $K$ contains a primitive $p$-th root of unity and
char$(\widehat K) = p$, this has been proved in \cite{BH}, Sect.
4, and in \cite{Ch1}, Sect. 2):
\par
\medskip
\begin{prop}
\label{prop2.2} Let $(K, v)$ be an {\rm HDV}-field with {\rm
char}$(\widehat K) = p > 0$. Then:
\par
{\rm (a)} {\rm Brd}$_{p}(K)$ is infinite if and only if $\widehat
K/\widehat K ^{p}$ is an infinite extension;
\par
{\rm (b)} {\rm Brd}$_{p}(K) \ge n$ if $[\widehat K\colon \widehat
K ^{p}] = p ^{n}$, for some $n \in \mathbb N$.
\end{prop}
\par
\medskip
When char$(\widehat K) = p > 0$, the inequality $[\widehat K\colon
\widehat K ^{p}] \le p$ holds if and only if $\widehat K$ is an
almost perfect field, i.e. its finite extensions are simple;
perfect fields of any characteristic satisfy the condition on the
right side and, in this sense, are also almost perfect (see
\cite{L}, Ch. V, Theorem~4.6 and Corollary~6.10). This allows us
to state the main results of the present paper as follows:
\par
\medskip
\begin{theo}
\label{theo2.3} Let $(K, v)$ be an {\rm HDV}-field with {\rm
char}$(\widehat K) = p > 0$. Then {\rm Brd}$_{p}(K) \le 1$ if and
only if $\widehat K$ is $p$-quasilocal and almost perfect; the
equality {\rm Brd}$_{p}(K) = 0$ holds if and only if $\widehat K$
is perfect and $\widehat K(p) = \widehat K$.
\end{theo}
\par
\medskip
Theorem \ref{theo2.3} yields Brd$_{p}(K) = 1$ in case $\widehat K
_{\rm sep} = \widehat K$ and $[\widehat K\colon \widehat K ^{p}] =
p$. This result is contained in \cite{Ya}, Proposition~2.1 (see
also \cite{BH}, Proposition~4.5), and it is used for proving the
stated theorem in general.
\par
\medskip
\begin{coro}
\label{coro2.4} Assuming that $(K, v)$ and $p$ satisfy the
conditions of Proposition \ref{prop2.2}, let $[\widehat K\colon
\widehat K ^{p}] = p$. Then {\rm Brd}$_{p}(K)  = 2$ unless
$\widehat K$ is $p$-quasilocal.
\end{coro}

\medskip
Corollary \ref{coro2.4} follows from Theorem \ref{theo2.3} and
\cite{PS}, Theorem~2. Theorem \ref{theo2.3} and this corollary
determine Brd$_{p}(K)$ when $(K, v)$ is HDV with char$(\widehat K)
= p$ and $\widehat K$ almost perfect. At the same time,
Proposition \ref{prop2.1} and Theorem \ref{theo2.3} yield the
following characterization of stable HDV-fields:
\par
\medskip
\begin{coro}
\label{coro2.5} Let $(K, v)$ be an {\rm HDV}-field. Then $K$ is
stable if and only if $\widehat K$ is almost perfect, stable and
$p$-quasilocal, for every $p \in \mathbb{P}$.
\end{coro}
\par
\medskip
It is presently unknown whether a field $E$ is stable, under the
condition that it is $p$-quasilocal, for every $p \in \mathbb{P}$.
In view of the above-noted facts on the Brauer $p$-dimensions of
$p$-quasilocal fields, the stability of $E$ will be proved if the
following open problem has an affirmative solution:
\par
\medskip
{\bf Problem 2.6.} Let $F$ be a field not containing a primitive
$p$-th root of unity, for some $p \in \mathbb{P}$ different from
char$(F)$. Find whether Brd$_{p}(F) = 0$ in case $F(p) = F$.
\par
\medskip
The basic notation, terminology and conventions kept in this paper
are standard and essentially the same as in \cite{L}, \cite{P},
\cite{JW} and \cite{Ch3}. We refer the reader to \cite{P}, for the
definition of a cyclic algebra over an arbitrary field; the
notions of an inertial algebra, an inertial lift, and a nicely
semiramified (briefly, NSR) algebra over a Henselian field are
defined in \cite{JW}. As in \cite{Ch7}, we suppose that for any
discrete valued field $(K, v)$, $v(K)$ is chosen to be a subgroup
of the additive group $\mathbb{Q}$ of rational numbers.
Throughout, Brauer groups and value groups are written additively,
Galois groups are viewed as profinite with respect to the Krull
topology, and by a profinite group homomorphism, we mean a
continuous one. Given a field $E$, $E ^{\ast }$ denotes its
multiplicative group, $E ^{\ast n} = \{a ^{n}\colon \ a \in E
^{\ast }\}$, for each $n \in \mathbb N$, and $\mathcal{G}_{E} =
\mathcal{G}(E _{\rm sep}/E)$ is the absolute Galois group of $E$.
For any $p \in \mathbb P$, we denote by $_{p}{\rm Br}(E)$ the
group $\{b _{p} \in {\rm Br}(E)\colon \ pb _{p} = 0\}$ and by $r
_{p}(E)$ the rank of $\mathcal{G}(E(p)/E)$ as a pro-$p$-group,
i.e. the cardinality of any minimal system of generators of
$\mathcal{G}(E(p)/E)$ as a topological group (we put $r _{p}(E) =
0$ if $E(p) = E$). As usual, Br$(E ^{\prime }/E)$ stands for the
relative Brauer group of any field extension $E ^{\prime }/E$
(defined to be the kernel of the scalar extension map $\pi _{E\to
E'}$ of Br$(E)$ into Br$(E ^{\prime })$); also, $E ^{\prime }$ is
called a splitting field of every $A \in s(E)$ with $[A] \in {\rm
Br}(E ^{\prime }/E)$. We write I$(E ^{\prime }/E)$ for the set of
intermediate fields of $E ^{\prime }/E$; when $E ^{\prime }/E$ is
separable of finite degree $[E ^{\prime }\colon E]$, $N(E ^{\prime
}/E)$ denotes the norm group of $E ^{\prime }/E$.
\par
\medskip
Here is an overview of the paper: Section 3 includes preliminaries
on Henselian fields used in the sequel. Theorem \ref{theo2.3} is
proved in Section 4. Absolutely stable HDV-fields are
characterized in Section 5, where some special fields of this kind
are also presented. Specifically, we show that an $m$-dimensional
local field $K _{m}$ (i.e. a complete $m$-discretely valued field
(see \cite{TY}, \cite{F1}, \cite{Zh}) with a quasifinite $m$-th
residue field) is absolutely stable, if $m \le 2$, and $K _{m}$ is
not stable, otherwise. When char$(K _{m}) > 0$, this is contained
in \cite{Ch1}, Corollaries~4.5, 4.6.
\par
\medskip

\section{Preliminaries}
\label{Sec:3}
\medskip
Let $K$ be a field with a nontrivial valuation $v$, $O _{v}(K) =
\{a \in K\colon \ v(a) \ge 0\}$ the valuation ring of $(K, v)$, $M
_{v}(K) = \{\mu \in K\colon \ v(\mu ) > 0\}$ the maximal ideal of
$O _{v}(K)$, $O _{v}(K) ^{\ast } = \{u \in K\colon \ v(u) = 0\}$
the group of units of $O _{v}(K)$, $v(K)$ and $\widehat K = O
_{v}(K)/M _{v}(K)$ the value group and the residue field of $(K,
v)$, respectively. For each $\gamma \in v(K)$, $\gamma \ge 0$,
$\nabla _{\gamma }(K)$ denotes the set $\{\lambda \in K\colon \
v(\lambda - 1) > \gamma \}$. Note that $v$ is Henselian if the
following conditions hold: $K$ is complete relative to the
topology of $v$; $v(K)$ is an Archimedean group, i.e. it embeds as
an ordered subgroup in the additive group $\mathbb R$ of real
numbers (cf. \cite{L}, Ch. XII, and see Lemma \ref{lemm3.4} below).
In order that $v$ be Henselian, it is necessary and sufficient
that any of the following two equivalent conditions holds (cf.
\cite{E3}, Sect. 18.1, and \cite{L}, Ch. XII, Sect. 4):
\par
\medskip
(3.1) (a) Given a polynomial $f(X) \in O _{v}(K) [X]$ and an
element $a \in O _{v}(K)$, such that $2v(f ^{\prime }(a)) <
v(f(a))$, where $f ^{\prime }$ is the formal derivative of $f$,
there is a zero $c \in O _{v}(K)$ of $f$ satisfying the equality
$v(c - a) = v(f(a)/f ^{\prime }(a))$;
\par
(b) For each normal extension $\Omega /K$, $v ^{\prime }(\tau (\mu
)) = v ^{\prime }(\mu )$ whenever  $\mu \in \Omega $, $v ^{\prime
}$ is a valuation of $\Omega $ extending $v$, and $\tau $ is a
$K$-automorphism of $\Omega $.
\par
\medskip
When $v$ is Henselian, so is $v _{L}$, for every algebraic field
extension $L/K$. In this case, we denote by $\widehat L$ the
residue field of $(L, v _{L})$, and put $O _{v}(L) = O _{v
_{L}}(L)$, $M _{v}(L) = M _{v_{L}}(L)$ and $v(L) = v _{L}(L)$.
Clearly, $\widehat L$ is an algebraic extension of $\widehat K$,
and $v(K)$ is an ordered subgroup of $v(L)$; the index of $v(K)$
in $v(L)$ is denoted by $e(L/K)$. By Ostrowski's theorem, if
$[L\colon K]$ is finite, then $[\widehat L\colon \widehat
K]e(L/K)$ divides $[L\colon K]$ and $[L\colon K][\widehat L\colon
\widehat K] ^{-1}e(L/K) ^{-1}$ has no divisor $p \in \mathbb P$,
$p \neq {\rm char}(\widehat K)$; in particular, $v _{L}$ can be
chosen so that $v(L)$ be an ordered subgroup of a fixed divisible
hull of $v(K)$. We say that $L/K$ is defectless, if $[L\colon K] =
[\widehat L\colon \widehat K]e(L/K)$. The defectlessness of $L/K$
is guaranteed, if char$(\widehat K) \nmid [L\colon K]$ as well as
in the following two cases:
\par
\medskip
(3.2) (a) If $(K, v)$ is HDV and $L/K$ is separable (see
\cite{E3}, Sect. 17.4);
\par
(b) When $(K, v)$ is a complete discrete valued field (cf.
\cite{L}, Ch. XII, Proposition~6.1).
\par
\medskip\noindent
Assume as above that $(K, v)$ is Henselian. We say that a finite
extension $R$ of $K$ is inertial, if $[R\colon K] = [\widehat
R\colon \widehat K]$ and $\widehat R/\widehat K$ is separable;
$R/K$ is called totally ramified if $e(R/K) = [R\colon K]$.
Inertial extensions are separable and have the following
properties (see \cite{TW}, Theorem~A.23):
\par
\medskip
\begin{lemm}
\label{lemm3.1} {\rm (a)} An inertial extension $R ^{\prime }/K$ is
Galois if and only if $\widehat R ^{\prime }/\widehat K$ is
Galois. When this holds, $\mathcal{G}(R ^{\prime }/K)$ and
$\mathcal{G}(\widehat R ^{\prime }/\widehat K)$ are canonically
isomorphic.
\par
{\rm (b)} The compositum $K _{\rm ur}$ of inertial extensions of
$K$ in $K _{\rm sep}$ is a Galois extension of $K$ with
$\mathcal{G}(K _{\rm ur}/K)$ isomorphic to $\mathcal{G}_{\widehat
K}$. Finite extensions of $K$ in $K _{\rm ur}$ are inertial, and
the natural map $I(K _{\rm ur}/K) \to I(\widehat K _{\rm
sep}/\widehat K)$ is bijective.
\par
{\rm (c)} The group $N(I/K)$ includes $\nabla _{0}(K)$, for every
inertial extension $I/K$.
\end{lemm}
\par
\medskip
The Henselian property of $(K, v)$ guarantees that $v$ extends on
each $D \in d(K)$ to a unique, up-to equivalence, valuation $v
_{D}$ (cf. \cite{JW}, pp. 131-132; \cite{TW}, Exercise~3.11). Put
$v(D) = v _{D}(D)$ and let $\widehat D$ be the residue division
ring of $(D, v _{D})$. It is known that $\widehat D$ is a division
$\widehat K$-algebra, $[\widehat D\colon \widehat K] < \infty $,
$v(D)$ is an ordered abelian group and $v(K)$ is an ordered
subgroup of $v(D)$ of finite index $e(D/K)$. Moreover, by
Ostrowski-Draxl's theorem \cite{Dr}, $[\widehat D\colon \widehat
K]e(D/K) \mid [D\colon K]$ and $[D\colon K][\widehat D\colon
\widehat K] ^{-1}e(D/K) ^{-1}$ has no prime divisor $p \neq {\rm
char}(\widehat K)$. In addition, Proposition~2.2 of \cite{TY},
states the following:
\par
\medskip
\begin{lemm}
\label{lemm3.2} Let $(K, v)$ be an {\rm HDV}-field and $D \in
d(K)$. Then $D/K$ is defectless, i.e. $[D\colon K] = [\widehat
D\colon \widehat K]e(D/K)$.
\end{lemm}
\par
\medskip\noindent
Next we  state results on inertial and central $K$-algebras
(contained in \cite{JW}, Theorem~2.8), which are used for studying
the sequence Brd$_{p}(K)$, $p \in \mathbb{P}$:
\par
\medskip
\begin{lemm}
\label{lemm3.3} Let $(K, v)$ be a Henselian field. Then the set
{\rm IBr}$(K) = \{[S ^{\prime }] \in {\rm Br}(K)\colon S ^{\prime
} \in d(K)$ is inertial over $K\}$ is a subgroup of {\rm Br}$(K)$,
and the natural map {\rm IBr}$(K) \to {\rm Br}(\widehat K)$ is an
index-preserving group isomorphism; {\rm Brd}$_{p}(\widehat K) \le
{\rm Brd}_{p}(K)$, for all $p \in \mathbb P$, and equality holds
if {\rm Brd}$_{p}(\widehat K) = \infty $.
\end{lemm}
\par
\medskip
The following lemma shows that a nontrivially valued field $(K,
v)$ with $v(K)$ Archimedean is Henselian if and only if $K$ does
not admit a separable proper extension in its completion $K _{v}$
with respect to the topology induced by $v$. This is a known
consequence of basic properties of valuation prolongations on
finite separable extensions (cf. \cite{L}, Ch. XII, Sects. 2, 3
and 6). When $(K, v)$ is Henselian, the lemma characterizes finite
extensions of $K _{v}$ in $K _{v,{\rm sep}}$. As it seems to be
difficult to find a standard reference to these results, we refer
the reader to \cite{Ch8}, Lemma~3.1, for a proof of the lemma.
\par
\medskip
\begin{lemm}
\label{lemm3.4} Assume that $(K, v)$ is a nontrivially valued field
with $v(K)$ Archimedean. Then $(K, v)$ is Henselian if and only if
$K$ is separably closed in $K _{v}$. When $(K, v)$ is Henselian,
the following conditions are fulfilled:
\par
{\rm (a)} Every $L \in {\rm Fe}(K _{v})$ is $K _{v}$-isomorphic to
$\widetilde L \otimes _{K} K _{v}$ and $\widetilde L _{v}$, where
$\widetilde L$ is the separable closure of $K$ in $L$. The
extension $L/K _{v}$ is Galois if and only if so is $\widetilde
L/K$; in case this holds, $\mathcal{G}(L/K _{v})$ and
$\mathcal{G}(\widetilde L/K)$ are isomorphic.
\par
{\rm (b)} $K _{\rm sep} \otimes _{K} K _{v}$ is a field and there
exist canonical isomorphisms $K _{\rm sep} \otimes _{K} K _{v}
\cong K _{v,{\rm sep}}$ and $\mathcal{G}_{K} \cong \mathcal{G}_{K
_{v}}$.
\end{lemm}
\par
\medskip
The proof of Theorem \ref{theo2.3} relies on the following
well-known results:
\par
\medskip
\begin{prop}
\label{prop3.5} Let $(K, v)$ be an {\rm HDV}-field, and $\bar v$
the valuation of $K _{v}$ continuously extending $v$. Then:
\par
{\rm (a)} The map $\pi _{K\to K _{v}}$ is an injective
homomorphism preserving Schur indices and exponents (cf.
\cite{Cohn}, Theorem~1), so {\rm Brd}$_{p'}(K) \le {\rm
Brd}_{p'}(K _{v})$, for every $p' \in \mathbb P$;
\par
{\rm (b)} The valued field $(K _{v}, \bar v)$ is an immediate
extension of $(K, v)$, i.e. it is a valued extension with
$\widehat K _{v} = \widehat K$ and $\bar v(K _{v}) = v(K)$ (cf.
\cite{E3}, Theorem~9.3.2);
\par
{\rm (c)} If char$(K) = p > 0$, $[\widehat K\colon \widehat K
^{p}] = p ^{n}$, for some integer $n \ge 0$, and $\mathcal{K}/K
_{v}$ is a finite extension, then $(\mathcal{K}, \bar v
_{\mathcal{K}})$ is a complete discrete valued field with
$[\mathcal{K}\colon \mathcal{K} ^{p}] = p ^{n+1}$ and $[\widehat
{\mathcal{K}}\colon \widehat {\mathcal{K}} ^{p}] = p ^{n}$ (apply
(3.2) (b), \cite{BH}, Lemma~2.12, and for the completeness of
$(\mathcal{K}, \bar v _{\mathcal{K}})$, see \cite{L}, Ch. XII,
Proposition~2.5).
\end{prop}
\par
\medskip
At the end of this Section, we prove the latter assertion of
Theorem \ref{theo2.3}. Let $(K, v)$ be an HDV-field with
char$(\widehat K) = p > 0$. Proposition \ref{prop2.2} allows to
consider only the case where $\widehat K$ is perfect. In this
case, Br$(\widehat K) _{p} = \{0\}$ (cf. \cite{A1}, Ch. VII,
Theorem~22), so Lemma \ref{lemm3.2} and \cite{JW}, Lemma~5.14,
imply that if $r _{p}(\widehat K) = 0$, then Br$(K) _{p} = \{0\}$,
i.e. Brd$_{p}(K) = 0$. Suppose that $r _{p}(\widehat K) > 0$. Then
there is an NSR-algebra $\Delta _{p} \in d(K)$ of degree $p$,
whence Brd$_{p}(K) \ge 1$. The inequality Brd$_{p}(K) \le 1$ is
known (see, e.g., \cite{PS}, Corollary~2.5); it is also a part of
the next lemma which we prove for convenience of the reader.
\par
\medskip
\begin{lemm}
\label{lemm3.6} Assume that $(K, v)$ is an {\rm HDV}-field, such
that $\widehat K$ is perfect, {\rm char}$(\widehat K) = p > 0$ and
$r _{p}(\widehat K)
> 0$. Then every $D _{p} \in d(K)$ of $p$-primary degree {\rm
deg}$(D _{p}) \neq 1$ is an {\rm NSR}-algebra over $K$ with {\rm
exp}$(D _{p}) = {\rm deg}(D _{p})$.
\end{lemm}
\par
\smallskip
\begin{proof}
Take any $K$-algebra $D _{p} \in d(K)$, $D _{p} \neq K$, and let
deg$(D _{p}) = p ^{n}$. As $\widehat K$ is perfect, whence
Br$(\widehat K) _{p} = \{0\}$, $\widehat D _{p}/\widehat K$ is a
separable field extension, so it is simple, which implies
$\widehat D _{p} = \widehat D _{p} ^{\prime }$ and $[D _{p}
^{\prime }\colon K] = [\widehat D _{p}\colon \widehat K]$ divides
$p ^{n}$, for some inertial extension $D _{p} ^{\prime }$ of $K$
included in $D _{p}$ as a (commutative) $K$-subalgebra (see
\cite{P},~Sect.~13.1). Since $v(K)$ is cyclic, it follows
similarly that $e(D _{p}/K)$ divides $p ^{n}$. Hence, by Lemma
\ref{lemm3.2}, $[\widehat D _{p}\colon \widehat K] = e(D _{p}/K) =
p ^{n}$. Therefore, by \cite{JW}, Corollary~6.10 (a result of
Platonov and Yanchevskij) and the cyclicity of the group $v(D
_{p})/v(K)$, $p ^{n} \mid {\rm exp}(D _{p})$, proving that exp$(D
_{p}) = p ^{n}$. Finally, one obtains from \cite{JW}, Lemma~5.14,
that $[D _{p}] = [N _{p}]$, where $N _{p} \in d(K)$ is NSR. By
Wedderburn's structure theorem (see, e.g., \cite{P}, Sect. 3.5),
this means that $D _{p} \cong N _{p}$ as $K$-algebras, so Lemma
\ref{lemm3.6} is proved.
\end{proof}

\par
\medskip
\section{Proof of the main theorem}
\par
\medskip
The purpose of this Section is to prove the former assertion of
Theorem \ref{theo2.3}. Let $(K, v)$ be an HDV-field with
char$(\widehat K) = p$. Using Propositions \ref{prop2.1}~(a),
\ref{prop2.2} and the latter assertion of Theorem \ref{theo2.3},
one sees that it suffices to deduce the inequality Brd$_{p}(K) \le
1$, provided that $[\widehat K\colon \widehat K ^{p}] = p$ and
$\widehat K$ is $p$-quasilocal. Proposition \ref{prop3.5} (c) and
the equality $[\widehat K\colon \widehat K ^{p}] = p$ imply finite
extensions of $\widehat K$ are almost perfect fields and $\widehat
K$ has a unique, up-to a $\widehat K$-isomorphism, purely
inseparable extension $\widehat K _{n}$ of degree $p ^{n}$, for
each $n \in \mathbb N$. Therefore, using \cite{A1}, Ch. VII,
Theorem~32, one obtains the following:
\par
\medskip
(4.1) (a) Each $\widetilde D _{p} \in d(\widehat K)$ with
$[\widetilde D _{p}] \in {\rm Br}(\widehat K) _{p}$ has a
splitting field that is a purely inseparable extension of
$\widehat K$ of degree equal to exp$(\widetilde D _{p})$; in
particular, deg$(\widetilde D _{p}) = {\rm exp}(\widetilde D
_{p})$, i.e. Brd$_{p}(\widehat K) \le 1$;
\par
(b) If $\widetilde \Delta _{p} \in d(\widehat K)$ and
exp$(\widetilde \Delta _{p}) = p$, then $\widetilde \Delta _{p}$
is a cyclic $\widehat K$-algebra (cf. \cite{P}, Sect. 15.5);
\par
(c) The inner group product $Y ^{\ast g}\nabla _{0}(Y)$ includes
$O _{v}(K) ^{\ast }$ in case $Y/K$ is a finite extension,
$[Y\colon K] = [\widehat Y\colon \widehat K] = g$ and $\widehat
Y/\widehat K$ is purely inseparable.
\par
\medskip
The proof of Theorem \ref{theo2.3} also relies on the following
lemma.
\par
\medskip
\begin{lemm}
\label{lemm4.1} Let $(K, v)$ be an {\rm HDV}-field with {\rm
char}$(\widehat K) = p$ and $[\widehat K\colon \widehat K ^{p}] =
p$, and let $Y/K$ be a field extension, such that $[Y\colon K] =
[\widehat Y\colon \widehat K] = p$. Suppose that $\widehat K$ is
$p$-quasilocal and $\widehat Y$ is normal over $\widehat K$. Then
{\rm Br}$(Y/K)$ includes the group $_{p}{\rm Br}(K) \cap {\rm
IBr}(K)$, and the homomorphism $\pi _{K\to Y}\colon {\rm Br}(K)
\to {\rm Br}(Y)$ maps {\rm Br}$(K) _{p} \cap {\rm IBr}(K)$
surjectively upon {\rm Br}$(Y) _{p} \cap {\rm IBr}(Y)$.
\end{lemm}
\par
\medskip
\begin{proof}
It follows from \cite{Ch3}, Theorem~4.1, and Albert-Hochschild's
theorem (cf. \cite{S1}, Ch. II, 2.2) that $\pi _{\widehat K\to
\widehat Y}$ maps Br$(\widehat K) _{p}$ surjectively upon
Br$(\widehat Y) _{p}$. At the same time, we have Br$(\widehat
Y/\widehat K) =$ $_{p}{\rm Br}(\widehat K)$, by \cite{Ch3},
Theorem~4.1, if $\widehat Y/\widehat K$ is separable, and by (4.1)
(a), when $\widehat Y/\widehat K$ is inseparable. Note further
that IBr$(Y)$ includes the image of IBr$(K)$ under $\pi _{K\to
Y}$, and the natural maps $r _{K\to \widehat K}\colon {\rm IBr}(K)
\to {\rm Br}(\widehat K)$ and $r _{Y\to \widehat Y}\colon {\rm
IBr}(Y) \to {\rm IBr}(\widehat Y)$, are index-preserving group
isomorphisms (see \cite{JW}, Theorems~5.6 and 2.8). Since $(\pi
_{\widehat K\to \widehat Y} \circ r _{K\to \widehat K}) ([D]) = (r
_{Y\to \widehat Y} \circ \pi _{K\to Y}) ([D])$ (in Br$(\widehat
Y)$) whenever $D \in d(K)$ is inertial over $K$, this enables one
to prove the latter part of the assertion of Lemma \ref{lemm4.1},
and the fact that ind$(D _{p} \otimes _{K} Y) = {\rm deg}(D
_{p})/p$, for each $D _{p} \in d(K)$ with $[D _{p}] \neq 0$ and
$[D _{p}] \in ({\rm Br}(K) _{p} \cap {\rm IBr}(K))$ (the stated
equality is also implied by (4.1), Lemma \ref{lemm3.1} and
\cite{P}, Sect. 15.1, Proposition~b). In view of the Corollary in
\cite{P}, Sect. 13.4, these results complete our proof.
\end{proof}
\par
\medskip
Next we show that Theorem \ref{theo2.3} will be proved, if we
deduce the equality deg$(\Delta ) = p$, assuming that $\Delta \in
d(K)$ and exp$(\Delta ) = p$. It follows from Lemma \ref{lemm3.2}
and \cite{JW}, Proposition~1.7, that each $D \in d(K)$ with
deg$(D) = p$ possesses a maximal subfield $Y$ satisfying the
conditions of Lemma \ref{lemm4.1}. Hence, $\widehat Y$ is
$p$-quasilocal (cf. \cite{Ch3}, Theorem~4.1 and Proposition~4.4),
which enables one to obtain from the claimed property of $\Delta
$, by the method of proving \cite{Ch4}, Lemma~4.1, that if $\Delta
_{n} \in d(K)$ and exp$(\Delta _{n}) = p ^{n}$, then $\Delta _{n}$
has a splitting field $Y _{n}$ with $[Y _{n}\colon K] = p ^{n}$,
$v(Y _{n}) = v(K)$ and $\widehat Y _{n} \in I(\widehat Y ^{\prime
}/\widehat K)$, where $\widehat Y ^{\prime }$ is a perfect closure
of $\widehat K(p)$. This result gives the desired reduction.
Since, by Merkur'ev's theorem \cite{M}, Sect. 4, Theorem~2, each
$\Delta \in d(K)$ with exp$(\Delta ) = p$ is Brauer equivalent to
a tensor product of degree $p$ algebras from $d(K)$, we need only
prove that if $D _{j} \in d(K)$ and deg$(D _{j}) = p$, $j = 1, 2$,
then $D _{1} \otimes _{K} D _{2} \notin d(K)$. This can be deduced
from the next lemma.
\par
\medskip
\begin{lemm}
\label{lemm4.2} Let $(K, v)$ be an {\rm HDV}-field with {\rm
char}$(\widehat K) = p$, $\widehat K$ $p$-quasilocal and
$[\widehat K\colon \widehat K ^{p}] = p$. Then {\rm exp}$(\Delta )
= p ^{2}$, for any $\Delta \in d(K)$ of degree $p ^{2}$.
\end{lemm}
\par
\medskip
\begin{proof}
Let $\Delta $ be a $K$-algebra satisfying the conditions of the
lemma. As $\widehat K$ is almost perfect, this implies $p ^{2}$ is
divisible by the dimension of any commutative $\widehat
K$-subalgebra of $\widehat \Delta $. At the same time, it follows
from Lemma \ref{lemm3.2} and the cyclicity of $v(\Delta )$ that
$e(\Delta /K) \mid p ^{2}$. Suppose first that $e(\Delta /K) = 1$.
Then $[\widehat \Delta \colon \widehat K] = p ^{4}$, by Lemma
\ref{lemm3.2}, so the observation on commutative $\widehat
K$-subalgebras of $\widehat \Delta $ indicates that $\widehat
\Delta \in d(\widehat K)$ and deg$(\widehat \Delta ) = p ^{2}$.
Applying \cite{A1}, Ch. VII, Theorem~28, and \cite{Ch3},
Theorem~3.1, one concludes that exp$(\widehat \Delta ) = p ^{2}$.
It is now easily obtained from \cite{JW}, Theorems~2.8 and 2.9,
that $\Delta /K$ is inertial and deg$(\Delta ) = {\rm exp}(\Delta
)$, as claimed by Lemma \ref{lemm4.2}.
\par
\medskip
Henceforth, we assume that $e(\Delta /K) \neq 1$. Our first goal
is to prove that:
\par
\medskip
(4.2) (a) If $U$ is a central $K$-subalgebra of $\Delta $ of
degree $p$, then $U$ is neither an inertial nor an NSR-algebra
over $K$;
\par
(b) If $e(\Delta /K) = p$, then totally ramified extensions of $K$
of degree $p$ are not embeddable in $\Delta $ as $K$-subalgebras.
\par
\medskip\noindent The proof of (4.2) (a) relies on the Double
Centralizer Theorem (see \cite{P}, Sect. 12.7), which implies
$\Delta $ is $K$-isomorphic to $U \otimes _{K} U ^{\prime }$, for
some $U ^{\prime } \in d(K)$ with deg$(U ^{\prime }) = p$. Suppose
for a moment that $U/K$ is inertial. Applying (3.2) (a), Lemma
\ref{lemm3.2} and \cite{JW}, Theorem~2.8 and Proposition~1.7, one
concludes that $e(U ^{\prime }/K) = p$, $\widehat U ^{\prime
}/\widehat K$ is a normal field extension of degree $p$, and $U
^{\prime }$ contains as a $K$-subalgebra an extension $Y$ of $K$
with $\widehat Y = \widehat U ^{\prime }$. Therefore, by Lemma
\ref{lemm4.1}, $Y$ is embeddable in $U$ as a $K$-subalgebra, which
means that $U \otimes _{K} Y \notin d(Y)$. As $\Delta \in d(K)$
and $U \otimes _{K} Y$ is a $K$-subalgebra of $\Delta $, this is a
contradiction proving that $U/K$ cannot be inertial. The
non-existence of an NSR-subalgebra of $\Delta $ of degree $p$ is
merely a consequence of (4.2) (b).
\par
We turn to the proof of (4.2) (b), so we assume that $e(\Delta /K)
= p$. Suppose that our assertion is false, i.e. $\Delta $ contains
as a $K$-subalgebra a totally ramified extension $T$ of $K$ of
degree $p$, and let $W ^{\prime }$ be the centralizer of $T$ in
$\Delta $. It is clear from the Double Centralizer Theorem that $W
^{\prime } \in d(T)$ and deg$(W ^{\prime }) = p$, and it follows
from Lemma \ref{lemm3.2} and the assumptions on $\Delta /K$ and
$T/K$ that $[\widehat W ^{\prime }\colon \widehat T] = p ^{2}$. As
$\widehat T = \widehat K$ is almost perfect, each commutative
$\widehat T$-subalgebra $\widehat \Theta ^{\prime }$ of $\widehat
W ^{\prime }$ embeds as a $\widehat T$-subalgebra in $\widehat
\Theta $, for some commutative $T$-subalgebra $\Theta $ of $W
^{\prime }$, so the noted facts show (similarly to the proof of
the equality $[\widehat D _{p}\colon \widehat T] = {\rm deg}(D
_{p})$, in the setting of Lemma \ref{lemm3.6}) that $[\widehat
\Theta ^{\prime }\colon \widehat T]$ equals $1$ or $p$. Thus they
prove that $\widehat W ^{\prime } \in d(\widehat T)$ and
deg$(\widehat W ^{\prime }) = p$. Taking again into account that
$\widehat T = \widehat K$, and using \cite{JW}, Theorem~2.8, one
concludes that $W ^{\prime } \cong W \otimes _{K} T$ as a
$T$-algebra, where $W \in d(K)$ is an inertial lift of $\widehat W
^{\prime }$ over $K$. This leads to the conclusion that $W$ is
embeddable in $\Delta $ as a $K$-subalgebra, which contradicts the
non-existence of inertial central $K$-subalgebras of $\Delta $ of
degree $p$. The obtained contradiction proves (4.2) (b) and
completes the proof of (4.2) (a).
\par
\vskip0.2truecm We continue with the proof of Lemma \ref{lemm4.2}
in the case of $e(\Delta /K) = p$. Clearly, Lemma \ref{lemm3.2}
yields $[\widehat \Delta \colon \widehat K] = p ^{3}$, so the
assumption that $[\widehat K\colon \widehat K ^{p}] = p$ implies
$\widehat \Delta $ is noncommutative. This means that $[\widehat
\Delta \colon Z(\widehat \Delta )] = p ^{2}$ and $[Z(\widehat
\Delta )\colon \widehat K] = p$, where $Z(\widehat \Delta )$ is
the centre of $\widehat \Delta $. First we prove that exp$(\Delta
) = p ^{2}$, under the extra hypothesis that $\Delta $ possesses a
$K$-subalgebra $\Delta _{0}$, such that $[\Delta _{0}\colon K] = p
^{3}$ and $\widehat \Delta _{0}$ is $\widehat K$-isomorphic to
$\widehat \Delta $; by \cite{JW}, Theorem~2.9, this holds in the
special case where $Z(\widehat K)$ is a separable extension of
$\widehat K$. It follows from \cite{JW}, Proposition~1.7, our
extra hypothesis and the cyclicity of $v(K)$ that $Z(\widehat
\Delta )/\widehat K$ is a normal extension of degree $p$. Hence,
by Lemma \ref{lemm4.1}, $[\Delta _{0}] = [D \otimes _{K} Z(\Delta
_{0})]$ (in Br$(Z(\Delta _{0}))$), for some $D \in d(K)$ inertial
over $K$. The obtained result shows that $[\Delta \otimes _{K} D
^{\rm op}] \in {\rm Br}(Z(\Delta _{0})/K)$, where $D ^{\rm op}$ is
the $K$-algebra opposite to $D$. This requires that exp$(\Delta
\otimes _{K} D ^{\rm op}) \mid p$. Since deg$(D) = {\rm exp}(D) =
p ^{2}$, it follows now that exp$(\Delta ) = p ^{2}$, as claimed.
\par
\vskip0.1truecm We are now prepared to consider the case of
$e(\Delta /K) = p$ in general. The preceding part of our proof
allows us to assume that $Z(\widehat \Delta )$ is a purely
inseparable extension of $\widehat K$. Note also that $[Z(\widehat
\Delta )\colon \widehat K] = p$, and it follows from \cite{Ch3},
Theorem~3.1, and \cite{A1}, Ch. VII, Theorem~28, that $\widehat
\Delta $ is a cyclic $Z(\widehat \Delta )$-algebra of degree $p$.
Therefore, there exists $\eta \in \Delta $, which generates an
inertial cyclic extension of $K$ of degree $p$. Hence, by the
Skolem-Noether theorem (cf. \cite{P}, Sect. 12.6), there is $\xi
\in \Delta ^{\ast }$, such that $\xi \eta ^{\prime }\xi ^{-1} =
\varphi (\eta ^{\prime })$, for every $\eta ^{\prime } \in K(\eta
)$, where $\varphi $ is a generator of $\mathcal{G}(K(\eta )/K)$.
Denote by $B$ the $K$-subalgebra of $\Delta $ generated by $\eta $
and $\xi $. It is easy to see that $K(\xi ^{p}) = Z(B)$, deg$(B) =
p$ and $B$ is either an inertial or an NSR-algebra over $K(\xi
^{p})$. In view of (4.2) (a), this means that $\xi ^{p} \notin K$
which gives $[K(\xi ^{p})\colon K] = p$, and combined with (4.2)
(b), proves that $v(K(\xi ^{p})) = v(K)$. In other words, $K(\xi
^{p}) ^{\ast } = O _{v}(K(\xi ^{p})) ^{\ast }.K ^{\ast }$. As
$e(\Delta /K) = p$, the obtained properties of $B$ and $K(\xi
^{p})$ indicate that if $B/K(\xi ^{p})$ is inertial (equivalently,
if $v _{\Delta }(\xi ) \in v(K)$, see \cite{JW}, Theorem~5.6 (a)),
then $\widehat B \cong \widehat \Delta $ over $\widehat K$. This
means that $\Delta /K$ is subject to the extra hypothesis, which
yields exp$(\Delta ) = p ^{2}$. When $B/K(\xi ^{p})$ is NSR, these
properties imply with (4.2) (b) and \cite{P}, Sect. 15.1,
Proposition~b, the existence of an algebra $\Sigma \in d(K)$
satisfying the following conditions:
\par
\medskip
(4.3) (a) $\Sigma $ is isomorphic to the cyclic $K$-algebra
$(K(\eta )/K, \varphi , \pi ^{\prime })$, for some $\pi ^{\prime }
\in K ^{\ast }$; $\Sigma /K$ is NSR, whence $\Sigma $ does not
embed in $\Delta $ as a $K$-subalgebra;
\par
(b) ind$(\Delta \otimes _{K} \Sigma ) = p ^{2}$ (see also
\cite{P}, Sect. 13.4, and \cite{Ch5}, (1.1)(b)), the underlying
division $K$-algebra $\Delta ^{\prime }$ of $\Delta \otimes _{K}
\Sigma $ has a $K$-subalgebra $Z ^{\prime }$ isomorphic to $Z(B)$,
and the centralizer $C _{\Delta '}(Z ^{\prime }) := C$ is an
inertial $Z ^{\prime }$-algebra.
\par
\medskip\noindent
Note here that $[\Delta ^{\prime }] \in {\rm Br}(K(\xi ^{p}, \eta
)/K)$. Using (3.2) (a), (4.3) and Lemma \ref{lemm3.2} (and also,
the Double Centralizer Theorem), one concludes that $[C\colon K] =
p ^{3}$ and either $\Delta ^{\prime }/K$ is inertial or $e(\Delta
^{\prime }/K) = p$ and $\widehat C \cong \widehat \Delta ^{\prime
}$ as a $\widehat K$-algebra. As shown above, this requires that
exp$(\Delta ^{\prime }) = p ^{2}$. In view of (4.3) (b) and the
equality deg$(\Sigma ) = {\rm exp}(\Sigma ) = p$, it thereby
proves that exp$(\Delta ) = p ^{2}$ as well.
\par
\medskip
It remains to consider the case where $e(\Delta /K) = p ^{2}$. We
first show that one may assume without loss of generality that
Brd$_{p}(\widehat K) = 0$. It follows from (4.2) (a), Lemma
\ref{lemm3.2} and the equality $e(\Delta /K) = p ^{2}$ that
$\widehat \Delta /\widehat K$ is a field extension of degree $p
^{2}$. Using \cite{JW}, Theorem~3.1, one obtains that $\Delta
\otimes _{K} U \in d(U)$, $v(\Delta \otimes _{K} U) = v(\Delta )$
and $e((\Delta \otimes _{K} U)/U) = p ^{2}$, provided $U$ is an
extension of $K$ in $K(p) \cap K _{\rm ur}$, such that no proper
extension of $\widehat K$ in $\widehat U$ is embeddable in
$\widehat \Delta $ as a $\widehat K$-subalgebra. Note also that
$\widehat \Delta \otimes _{\widehat K} \widehat U$ is $\widehat
U$-isomorphic to the residue field of $\Delta \otimes _{K} U$,
which enables one to prove (by applying Galois theory and Zorn's
lemma) that $U$ can be chosen so that $r _{p}(\widehat U) \le 1$.
Then, by \cite{J}, Proposition~4.4.8, Br$(\widehat U) _{p} =
\{0\}$, which leads to the desired reduction.
\par
\vskip0.1truecm We suppose further that Brd$_{p}(\widehat K) = 0$
and prove the following assertion:
\par
\medskip
(4.4) If $\Delta $ possesses a $K$-subalgebra $Z$, such that
$[Z\colon K] = [\widehat Z\colon \widehat K] = p$ and $\widehat Z$
is purely inseparable over $\widehat K$, then $\widehat \Delta
/\widehat K$ is purely inseparable.
\par
\medskip\noindent
Assuming the opposite and using (3.2) (a) and Lemma \ref{lemm3.2},
one obtains that $Z$ has an inertial extension $M$ which is a
maximal subfield of $\Delta $. As $v$ is Henselian, the
assumptions on $Z$ and $M$ ensure that $M = LZ$, for some inertial
extension $L$ of $K$ in $M$ of degree $p$. Note further that
$[M\colon K]$, $[\widehat M\colon \widehat K]$ and $[\widehat
\Delta\colon \widehat K]$ are equal to $p ^{2}$, which means that
$\widehat M = \widehat \Delta $. This enables one to deduce from
Lemma \ref{lemm3.1} and \cite{JW}, Proposition~1.7, that $L/K$ is
a cyclic extension. At the same time, the equality
Brd$_{p}(\widehat K) = 0$ and Albert-Hochschild's theorem, applied
to the extension $\widehat Z/\widehat K$, indicate that
Brd$_{p}(\widehat Z) = 0$. Therefore, the group $N(M/Z)$ includes
$O _{v}(Z) ^{\ast }$ (cf. \cite{P}, Sect. 15.1, Proposition~b),
which allows to obtain from the Skolem-Noether theorem and the
Double Centralizer Theorem that there is a $Z$-isomorphism $C
_{\Delta }(Z) \cong (M/Z, \psi ^{\prime }, \gamma )$, for some
$\gamma \in K ^{\ast }$ and a generator $\psi ^{\prime }$ of
$\mathcal{G}(M/Z)$. This implies $\Delta \cong D _{1} \otimes _{K}
D _{2}$ as $K$-algebras, where $D _{1} = (L/K, \psi , \gamma )$,
$\psi $ is the $K$-automorphism of $L$ induced by $\psi ^{\prime
}$, $D _{2} \in d(K)$ and $[D _{2}] \in {\rm Br}(Z/K)$. As
Brd$_{p}(\widehat K) = 0$, $\widehat K$ is almost perfect and
deg$(D _{2}) = p$, one obtains further that $D _{2}$ has a
$K$-subalgebra $T$ that is a totally ramified extension of $K$ of
degree $p$. It is easy to see that the $T$-algebra $L \otimes _{K}
T$ is a field. More precisely, $(L \otimes _{K} T)/T$ is an
inertial and cyclic extension of degree $p$, which allows to
deduce consecutively that the norm group $N((L \otimes _{K} T)/T)$
includes $O _{v}(T) ^{\ast }$ and $K ^{\ast }$. Observing also
that $D _{1} \otimes _{K} T$ is $T$-isomorphic to $((L \otimes
_{K} T)/T, \psi _{T}, \gamma )$, where $\psi _{T}$ is the
$T$-isomorphism of $L \otimes _{K} T$ extending $\psi $, one
obtains from \cite{P}, Sect. 15.1, Proposition~b, that $D _{1}
\otimes _{K} T \in s(T) \setminus d(T)$, whence, $D _{1} \otimes
_{K} T$ contains zero-divisors. As $D _{1} \otimes _{K} T$ is a
$K$-subalgebra of $D _{1} \otimes _{K} D _{2}$ and $D _{1} \otimes
_{K} D _{2} \cong \Delta \in d(K)$, this is a contradiction
proving (4.4).
\par
We are now in a position to prove Lemma \ref{lemm4.2}. If $\widehat
\Delta /\widehat K$ is a purely inseparable field extension, then
it follows from Proposition \ref{prop3.5} and \cite{Ya},
Proposition~2.1, that exp$(\Delta ) = p ^{2}$. Suppose finally
that $\widehat \Delta $ is a field and $\widehat \Delta /\widehat
K$ is not purely inseparable. In view of \cite{JW},
Proposition~1.7 and Theorem~2.9, this ensures the existence of an
inertial cyclic extension $\Lambda $ of $K$ of degree $p$, which
embeds in $\Delta $ as a $K$-subalgebra. Our goal is to show that
there is an infinite extension $W$ of $K$ in an algebraic closure
$\overline K$, satisfying the following:
\par
\vskip0.2truecm (4.5) $v(W) = v(K)$, $\widehat W$ is purely
inseparable over $\widehat K$ and $\Delta \otimes _{K} W \in
d(W)$.
\par
\vskip0.2truecm\noindent Note that (4.5) implies exp$(\Delta ) = p
^{2}$. Indeed, $[\widehat K\colon \widehat K ^{p}] = p$, so it
follows from (3.2) (a) and (4.5) that $\widehat W$ is perfect;
hence, by Lemma \ref{lemm3.6}, $(\Delta \otimes _{K} W)/W$ is NSR
and exp$(\Delta \otimes _{K} W) = {\rm deg}(\Delta \otimes _{K} W)
= p ^{2}$. Since exp$(\Delta \otimes _{K} W) \mid {\rm exp}(\Delta
)$ and exp$(\Delta ) \mid {\rm deg}(\Delta ) = p ^{2}$, this gives
exp$(\Delta ) = p ^{2}$, as required.
\par
Finally, we prove (4.5). Fix an element $a _{0} \in O _{v}(K)
^{\ast }$ so that $\hat a _{0} \notin \widehat K ^{p}$, take a
system $a _{n} \in \overline K$, $n \in \mathbb N$, satisfying $a
_{n} ^{p} = a _{n-1}$, for each $n$, and let $W$ be the union of
the fields $W _{n} = K(a _{n})$, $n \in \mathbb N$. It is easily
verified that $[W _{n}\colon K] = [\widehat W _{n}\colon \widehat
K] = p ^{n}$ and $\widehat W _{n}/\widehat K$ is purely
inseparable, for every $n \in \mathbb N$, so it follows from (3.2)
(a), the equality $[\widehat K\colon \widehat K ^{p}] = p$ and the
inclusions $W _{n} \subset W _{n+1}$, $n \in \mathbb N$, that $W$
is a field, $v(W) = v(K)$ and $\widehat W$ is a perfect closure of
$\widehat K$. Arguing by induction on $n$, taking into account
that $\Delta \otimes _{K} W _{n+1}$ and $(\Delta \otimes _{K} W
_{n}) \otimes _{W _{n}} W _{n+1}$ are isomorphic as $W
_{n+1}$-algebras, and using (4.4), the noted properties of $W
_{n}$, and the behaviour of Schur indices under scalar extensions
of finite degrees (cf. \cite{P}, Sect. 13.4), one obtains that,
for each $n \in \mathbb N$, $\Delta \otimes _{K} W _{n} \in d(W
_{n})$, and $\Lambda \otimes _{K} W _{n}$ is an inertial cyclic
extension of $W _{n}$ of degree $p$, embeddable in $\Delta \otimes
_{K} W _{n}$ as a $W _{n}$-subalgebra. Therefore, $\Delta \otimes
_{K} W \in d(W)$, so (4.5), Lemma \ref{lemm4.2} and Theorem
\ref{theo2.3} are proved.
\end{proof}
\par
\smallskip
\section{\bf Absolutely stable HDV-fields}
\par
\medskip The first result of this Section states the following:
\par
\medskip
\begin{coro}
\label{coro5.1} Let $(K, v)$ be an {\rm HDV}-field. Then $K$ is
absolutely stable if and only if $\widehat K$ is quasilocal and
almost perfect; for instance, this holds when $\widehat K$ is a
complete discrete valued field with a quasifinite residue field.
\end{coro}
\par
\smallskip
\begin{proof}
Our former conclusion follows from Corollary \ref{coro2.4}, the
absolute stability of quasilocal fields, and the fact that the
class of almost perfect fields is closed under the formation of
finite extensions. Since complete discrete valued fields with
quasifinite residue fields are quasilocal and almost perfect (see
\cite{Se}, Ch. XIII, Sect. 3, and \cite{E3}, Theorem~12.2.3), our
latter conclusion is an immediate consequence of the former one.
\end{proof}
\par
\medskip
The latter part of Corollary \ref{coro5.1} can be restated by
saying that $2$-dimensional local fields $K _{2}$ with quasifinite
second residue fields $K _{0}$ are absolutely stable. This result
can be specified as follows:
\par
\medskip
\begin{prop}
\label{prop5.2} An $m$-dimensional local field $K _{m}$ with a
quasifinite $m$-th residue field $K _{0}$ is stable if and only if
$m \le 2$. When $m \le 2$, $K _{m}$ is absolutely stable and {\rm
Brd}$_{p'}(K _{m} ^{\prime }) = 1$, $p ^{\prime } \in \mathbb P$,
for every finite extension $K _{m} ^{\prime }/K _{m}$.
\end{prop}
\par
\medskip
\begin{proof}
It is known (cf. \cite{Se}, Ch. XIII, Sect. 3) that if $m = 1$,
then $K _{m}$ is a quasilocal field with Br$(K _{m})$ isomorphic
to the quotient group $\mathbb Q/\mathbb Z$ of the additive group
of rational numbers by the subgroup of integers. This implies $K
_{m}$ is absolutely stable and Brd$_{p'}(K _{m} ^{\prime }) = 1$,
for all $p' \in \mathbb P$ and every finite extension $K _{m}
^{\prime }/K _{m}$, as claimed. We assume further that $m \ge
2$. Then $K _{m}$ is complete with respect to a discrete valuation $w
_{m}$ whose residue field $K _{m-1}$ is an $(m - 1)$-dimensional
local field with last residue field isomorphic to $K _{0}$. Therefore,
$(K _{m}, w _{m})$ is HDV, and by Lemma \ref{lemm3.3}, Brd$_{p'}(K
_{m-1}) \le {\rm Brd}_{p}(K _{m})$, for each $p ^{\prime } \in
\mathbb P$. Suppose now that $m = 2$. As noted above, $K _{1}$ is
quasilocal with Br$(K _{1}) \cong \mathbb Q/\mathbb Z$; in
addition, if char$(K _{1}) = {\rm char}(K _{0})$, then $K _{1}$ is
isomorphic to the formal power series field $K _{0}((X _{1}))$
(see \cite{E3}, Theorem~12.2.3), which is almost perfect. Hence,
by Corollary \ref{coro5.1}, $K _{2}$ is absolutely stable and
Brd$_{p'}(K _{u}) = 1$, $u = 1, 2$, $p' \in \mathbb P$. This
proves Proposition \ref{prop5.2}, for $m = 2$ (since finite
extensions of $K _{2}$ are $2$-dimensional local fields with
quasifinite 2nd residue fields). Next we prove that $r _{2}(K
_{1}) \ge 2$. Firstly, if char$(K _{0}) = 2$, then \cite{Ch8},
Lemma~2.2 and \cite{Ch5}, Lemma~4.2, show that $r _{2}(K _{1}) =
\infty $ unless char$(K _{1}) = 0$ and $K _{0}$ is finite.
Secondly, if char$(K _{1}) = 0$, char$(K _{0}) = 2$ and $K _{0}$
is finite, then it follows from Lemma \ref{lemm3.4} and \cite{S1},
Ch. II, Theorem~4, that $r _{2}(K _{1}) \ge 3$. When char$(K _{0})
\neq 2$, $K _{1} ^{\ast }/K _{1} ^{\ast 2}$ is a noncyclic group
of order $4$ (it is isomorphic to the direct sum $K _{0} ^{\ast
}/K _{0} ^{\ast 2} \oplus w _{1}(K _{1})/2w _{1}(K _{1})$, $w
_{1}$ being the discrete Henselian valuation of $K _{1}$ with
$\widehat K _{1} = K _{0}$), so it is clear from Kummer theory
that $r _{2}(K _{1}) = 2$. Lemma \ref{lemm3.1} and the inequality
$r _{2}(K _{1}) \ge 2$ imply there exist an algebra $\Delta _{2}
\in d(K _{2})$ and a field extension $L _{2}/K _{2}$, such that
deg$(\Delta _{2}) = [L _{2}\colon K _{2}] = 2$, $\Delta _{2}/K
_{2}$ is NSR, $L _{2}/K _{2}$ is inertial relative to $w _{2}$,
and $\Delta _{2} \otimes _{K} L _{2} \in d(L _{2})$. Thus it
follows that $K _{2}$ is not $2$-quasilocal. Using this result,
Proposition \ref{prop2.1} (a) and Lemma \ref{lemm3.3}, one obtains
that if $m \ge 3$, then Brd$_{2}(K _{j}) \ge 2$, $j = 3, \dots ,
m$, which completes our proof.
\end{proof}
\par
\medskip
In the setting of Proposition \ref{prop5.2}, $K _{m}$ is
$p$-quasilocal with Brd$_{p}(K _{m}) = 1$, provided that $p \in
\mathbb{P}$, $p \neq {\rm char}(K _{0})$ and $K _{0}$ does not
contain a primitive $p$-th root of unity (apply Proposition
\ref{prop2.1} (a) and \cite{Ch6}, Corollary~4.3, to the HDV-field
$(K _{m}, w _{m})$ with $\widehat K _{m} \cong K _{m-1}$, for $m
\ge 2$). When $K _{0}$ is finite, this means that Brd$_{p}(K _{m})
= 1$, for all $p \in \mathbb{P}$, with finitely many exceptions in
case $m \ge 3$, such as $p = 2$ and $p = {\rm char}(K _{0})$ (see
\cite{Ch6}, Proposition~4.4).
\par
\medskip
\begin{rema}
\label{rema5.3} Here are two special cases of Corollary
\ref{coro5.1} obtained in \cite{Ch1}:
\par
{\rm (i)} An {\rm HDV}-field $(K, v)$ with $\widehat K$ perfect is
absolutely stable if and only if $\widehat K$ is quasilocal
\cite{Ch1}, Corollary~4.6;
\par
{\rm (ii)} For any complete discrete valued field $(L, \omega )$
with a quasifinite residue field $\widehat L$ (specifically, for
any local field $L$), the formal power series field $L((T))$ is
absolutely stable \cite{Ch1}, Corollary~4.5 (ii). If char$(L) =
0$, this is also contained in \cite{Ch1}, Corollary~4.6. When
char$(L) = p > 0$, $L((T))$ is isomorphic to the iterated formal
power series field $\widehat L((Z))((T))$ (apply \cite{E3},
Theorem~12.2.3), so the inequality abrd$_{p}(L((T))) \le 1$, used
for proving \cite{Ch1}, Corollary~4.5 (ii), follows from
\cite{A1}, Ch. XI, Theorem~3, and results of Aravire, Jacob,
Merkurjev and Tignol (see \cite{AJ}, Sect. 3 and the Appendix).
\end{rema}
\par
\medskip
The concluding result of this paper is new if char$(K) \neq {\rm
char}(\widehat K)$ and $\widehat K$ is an imperfect field of type
$C _{1}$, in the sense of Lang and \cite{S1}, Ch. II. Under the
same hypotheses on $\widehat K$, if char$(K) = {\rm char}(\widehat
K)$, then the result is contained in \cite{Zheg}, Theorem~2, and
in case $\widehat K$ is perfect, it follows from \cite{Ch1},
Corollary~4.6.
\par
\medskip
\begin{coro}
\label{coro5.4} An {\rm HDV}-field $(K, v)$ is absolutely stable,
if $\widehat K$ has type $C _{1}$.
\end{coro}
\par
\smallskip
\begin{proof}
The field $\widehat K$ is almost perfect with abrd$_{p}(\widehat
K) = 0\colon p \in \mathbb P$ (cf. \cite{S1}, Ch. II, 3.2), so
$\widehat K$ is quasilocal, and by Corollary \ref{coro5.1}, $K$ is
absolutely stable.
\end{proof}
\par
\medskip

\emph{Acknowledgement}
This research has partially been supported by the Bulgarian
National Science Fund under Grant KP-06 N 32/1 of 07.12.2019.

\end{document}